# Compatibility of equations with truncated Newton's binomials.

## Anatoly A. Grinberg

## anatoly_gr@yahoo.com

### ABSTRACT

The resolvability of equations in integers containing truncated Newton's binomials, is determined by the divisibility of the binomial by the characteristic parameter of the equation, which most often is the binomial exponent. Two types of equations containing binomials from two and three integers are investigated. Conditions of resolvability of the equations are specified based on the characteristics of their parameters.

Equations containing truncated Newton's binomials , are solvable only at the certain congruence its terms. Criteria of congruence can be, as is usual, such properties as parity, co-primacy or divisibility by the characteristic parameters of a problem. Properties of divisibility of truncated binomials from two and three integers by their exponent n were described in article [1]. We use these results in the present work, assuming that n is a prime integer. In the first section we will describe equations containing a truncated binomial from the sum of two numbers. The second section is devoted to equations with binomials from three numbers. Truncated binomials are defined as Newton's series of exponent n from the sum of numbers $a, b, c$, from which the terms $a^n$, $b^n$, $c^n$ have been removed.

## I.Equations with truncated binomials of two integers.

Let us consider the equality

$$(a + b)^n = U(a, b) \qquad (I.1)$$

where $U(a, b)$- is the truncated binomial, which can be presented in

following three forms

$$U(a, b) = \sum_{\nu=1}^{n-1} \binom{n}{\nu} a^\nu b^{n-\nu} = -\sum_{\nu=1}^{n-1} \binom{n}{\nu} q^\nu (-a)^{n-\nu} = -\sum_{\nu=1}^{n-1} \binom{n}{\nu} q^\nu (-b)^{n-\nu} \qquad (I.2)$$

where $a$ and $b$ are integers, q = $a$ + b , $\binom{n}{\nu}$ are binomial coefficients proportional to a prime binomial exponent n. Therefore U (a, b) is proportional to n. In addition, U (a, b) is an even integer given any parity of integers $a$ and $b$. Without loss of generality it is possible to divide numbers $a$ and $b$ by their greatest common divider (GCD). As a result of such division, $a$ and $b$ are co-prime integers.

For compatibility with the right side of equation (I.1), the left side of the equation must also be even. Therefore, the parity of $a$ and $b$ must be the same. Because $a$ and $b$ are co-prime, they can only be odd.

Since 2 is a prime number, the divisibility of the left side of equation (I.1) by 2n, leads to the divisibility of $(a+b)$ by 2n. This brings us to the relation

$$q = (a + b) = 2\beta n \qquad (I.3)$$

where $\beta$ is parameter. Substituting (I.3) in (I.1) we obtain

$$(2\beta n)^n = U(a, b) \qquad (I.4)$$

The question of the possible implementation of equation (I.1) in integers, leads to the same question with respect to equation (I.4). While in equation (I.1) we could simply determine that both $a$ and b must be odd, equation (I.4) indicates a stronger criterion: its rght side must be divisible by $(2n)^n$ . This condition is necessary but not sufficient for Eq.(I.1) to exist in integers. If this division is possible, we must also show that there is no contradiction in the equation after division. If the division is not possible, then it is a clear indication of the impossibility of the Eq. (I.1) in integers.

Let us to formulate the above in a different order. If one needs to prove that the equation (I.4) cannot exist in integers, and it is known that U (a, b) is divisible by $(2n)^n$, , the contradiction must be sought in the equation remaining after the division.

On the other hand, the impossibility of existence of Eq.(I.4) in integers is not sufficient evidence of the indivisibility of the truncated binomial $U(a, b)$ by $(2n)^n$, [1].

Returning to the specifics of Eq.(I.1), it is easy to see that the difference between the left and right sides of the equation is given by $(a^n + b^n)$. Thus, it is always different from zero, except for the trivial case $a = -b$ . Therefore, the Eq.(I.1) does not exist in integers.

With regard to the divisibility of $U(a, b)$ by n, unfortunately, we do not know of an analysis of this property in the literature outside of its relation to the problem considered here. There are suppositions [1], based on the use of Fermat's little theorem, that the probability of their divisibility by $n^k$ is small for $k \geq 2$, but it does not solve the problem considered here.

─────────────────────────────

## II. Equation with a truncated binomial of three integers.

In this section we consider the more interesting case of the equation containing a truncated binomial of three integers $a$, $b$, c . This equation has the form

$$(a + b + c)^n = U(a, b, c) \qquad (II.1)$$

where

$$U(a, b, c) = U(a, b) + U((a + b), c) \qquad (II.2)$$

is a truncated three-integer binomial, which is, according to (II.2), comprised of the sum of two, two-integer truncated binomials. The divisibility of $U(a, b, c)$ by the binomial exponent was discussed in [1]. As in Section I, we are interested in the conditions under which equation (II.1) exists, at integers $a$, $b$, $c$ . Here we also will consider the exponent n to be a prime integer and suppose that the GCD has been removed from integers $a$, $b$, c .

Removing the GCD makes these numbers relatively co-prime, since the presence of a common factor in two of them, leads, in view of equation (II.1), to the presence of the same factor in the third. From the point of view of the divisibility of these numbers by exponent n, there are only two cases: the case where none of the numbers are divisible by n; and a case where one of them is divisible by a certain power of n ( the latter must be coordinated with the equation as well as with the dependence on n of other parameters). Since both binomials $U(a, b)$ and $U((a + b), c)$ are even, then from Eq.(II.1) and from the fact that $a$, $b$, $c$ are co-prime numbers, it follows that the two of them must be odd. We will assume that integer $a$ is even_. As in the case of Eq.(I.1), we obtain the equation

$$(a + b + c) = 2\beta n \qquad (II.3)$$

By substituting it in (II.1),we find

$$(2\beta n)^n = U(a, b) + U((a + b), c) \qquad (II.4)$$

The connection of the divisibility of the right side of equation (II.4) by $(2n)^n$ , with the existence in integers of equations like (II.1) or (II.4), was discussed in Section I. We consider now in detail the two above-mentioned cases of divisibility by n of numbers $a$, $b$, c .

## A. The case when numbers a, b, c are not divisible by n.

The question of the divisibility of truncated binomial $U(a, b)$ was discussed in [1] and in section I, where the low probability that it is divisible by $n^2$ was noted. Regarding binomial U((a + b), c), the sum of its arguments is, in accordance with (II.3), divisible by n. As shown in [1], in this case, U ((a + b), c) is divisible by $n^2$.

Thus, if $U(a, b)$ is divisible only by the first power of n, then the right-hand side of equation (II.4) is only divisible by n. On the other hand, if $U(a, b)$ is divisible by a power of n greater than 2, then the divisibility of the right side of (II.4) would be determined by the divisibility of $U((a + b), c)$. Uncertainty arises only when $U(a, b)$ is divisible by $n^2$. Thus, it is possible to formulate the following rule:

If none of the integers a, b, c is divisible by n and U(a, b) is divisible only by $n^k$, where , the equation (II. 1) cannot exist in integers

From this statement it follows that the question of the admissibility of equation (II.1) in integers, for certain values of exponent n can be easily solved by numerical methods, using modular arithmetic. In numerical calculations one only need clarify whether or not $U(a, b)$ is divisible by $n^2$. Below we will demonstrate this using cases with specific values of n

### 1. Equation (II.4) at n = 3.

$$(2\beta n)^n = n(a + b)(ab + 2\beta cn) \qquad (II.5)$$

According to Eq.(II.3), the sum $(a + b)$=q is not divisible by n, so $U(a, b)$ in (II.5) is divisible only by first power of n. Hence, equation does not exist in integers.

### 2. Equation (II.4) at n = 5.

$$(2\beta n)^n = \text{н}n(a + b)[ab (a^2 + ab + b^2) + 2 \beta cn((a + b)^2 + 2\beta cn)] \qquad (II.6)$$

Divisibility of the factor $(a^2 + ab + b^2)$ in (II.6) by n=5 is conveniently determined by congruence relations

$$(a^2 + ab + b^2) \equiv \Delta_a^2 + \Delta_a\Delta_b + \Delta_b^2 \pmod 5 \qquad (II.7)$$

where the remainders of integers $a$ and b, are indicated as $\Delta_a$, and $\Delta_b$ respectively. These vary from one to four. No combination of these remainders leads to a zero value or a multiple of n value of the right side of the (II.7), except when $\Delta_a + \Delta_b = n$. The latter corresponds to the multiplicity of $(a + b)$ to n, which is excluded from the conditions discussed here.

### 3. Equation (II.4) at n = 7.

$$(2\beta n)^n = n(a + b)[ab(a^2 + ab + b^2)^2 + 2 \beta cn((a + b)^2 + 2\beta cn)^2] \qquad (II.8)$$

This case is remarkable for the fact that it does not require any calculations. The problem of the existence of this equation is solved unambiguously for any divisibility of the factor $(a^2 + ab + b^2)$. Indeed, if this factor is not divisible by n, then $U(a, b)$ is proportional to the first power of n, and this is the divisibility of the right side of Eq.(II.8). If the same factor is divisible by n, then the divisibility of the right side of Eq.(II.8) is determined by the second term in its square brackets. It is interesting to compare this with another approach to the problem [2-4].

4. **Equation (II.4) at n = 11** .

The term U(a,b) of Eq.(II.4) at n=11 is given by

$$U(a, b) = nab\{ 5ab(a^7 + b^7) + 15a^2b^2(a^5 + b^5) +$$

$$30a^3b^3(a^3 + b^3) + 42a^4b^4(a + b) + a^9 + b^9 \} \qquad (II.9)$$

Its divisibility by $n^2$ was tested numerically using an modular arithmetic. The calculation showed that $U(a, b)$ is not divisible by $11^2$ for any numbers $a, b$ and therefore the terms of Eq.(II.4) are incompatible at n = 11.

## B. The case when one of the numbers $a$, $b$, $c$ is a multiple of n.

We assume that the numbers $c$ is a multiple of n raised to power $\rho_c$ and represent it in the form $c = c_0 n^{\rho_c}$ From (II.3) it follows that $q = (a + b)$ is also a multiple of some power of n. We denote this power as $\rho_q$. We also must attribute a dependence of parameter β on n, writing it in the form $\beta = \beta_0 n^{\rho_\beta}$. Thus, we have the following notation

$$c = c_0 n^{\rho_c}, \qquad q = q_0 \ n^{\rho_q}, \qquad \beta = \beta_0 n^{\rho_\beta} \qquad (II.10)$$

Inasmuch as we can solve the problem of the existence of Eq.(II,1) for at least one triplet of integers, we are free in the choice of the parameters $\rho_q$, $\rho_\beta$, $\rho_c$. We adjust them so that, if it is possible, they will provide the necessary divisibility of the right side of Eq.(II.1) by n. Substituting (II.10) into (II.3) we obtain

$$(q_0 n^{\rho_q} + c_0 n^{\rho_c}) = 2\beta_0 n \ n^{\rho_\beta} \qquad (II.11)$$

Let us assume that $\rho_c < \rho_q$. Since the first power of n is preserved when $\rho_c = 0$, from Eq.(II.11) follows the equality $\rho_c = \rho_\beta + 1$, which is valid only at $\rho_\beta > 0$.

On the other hand, from expansion of $U(a,b)$, as determined by Eq.(I.2), it can be seen that U(a, b) is proportional to $qn = q_0 n^{\rho_q+1}$. Similarly we find $U(a + b, c) \sim n^{\rho_q+1+\rho_c(n-1)}$.

From Eq. (II.4) it follows that $(2\beta_0 n^{\rho_\beta} n)^n \sim q_0 n^{\rho_q+1}$. Therefore $(\rho_\beta + 1)n = \rho_q + 1$. Again, this equality is preserved only if $\rho_c > 0$. Thus we arrive at the relations

$$\rho_\beta = \rho_c - 1 \ , \qquad \rho_q = n\rho_c - 1 \qquad \text{if } \rho_c > 0 \qquad\qquad (\text{II.12})$$

and the numbers c , q and β can be written as

$$c = c_0 n^{\rho_c} \ , \quad q = q_0\, n^{n\rho_c - 1} \ , \quad \beta = \beta_0 n^{\rho_c - 1} \qquad\qquad (\text{II.13})$$

From the above analysis it follows that, in the present case, we can always achieve the conditions under which the right-hand side of Eq.(II.4) is divisible by the necessary power of n. As discussed in section I, it does not mean that the Eq.(II.1) can exist in integers.

## III. Conclusion.

We examined a specific class of equations consisting of truncated Newtonian binomials of two and three numbers (trinomials). The main purpose was to establish the compatibility of terms of an equation for integer values of under-binomial numbers. The divisibility of the terms of the equations by the binominal exponent was choose as the main criterion of compatibility. The absence of necessary divisibility indicates incompatibility of the terms of the equation, while its presence leaves the question open. Limitation of results in this direction of inquiry is due to lack of theoretical development of the problem of truncated binomial divisibility. However, for certain values of parameters of the binomial, we have demonstrated a surprisingly simple solution to a problem that requires extensive analytical calculations when approached with conventional methods. Hopefully, this work will draw attention to solving the problem of the divisibility of truncated binomial U (a, b).

The author is grateful to B.M. Ashkinadze for useful discussions and valuable support. He also takes this opportunity to express his deep appreciation for I.D. Shkredov consulting in matters of number theory, and to B.B. Lurje for useful, critical remarks.